\documentclass[12pt]{amsart}

\usepackage{amsfonts,amsmath,latexsym,amssymb,verbatim,amsbsy,amsthm}

\usepackage[utf8]{inputenc}
\usepackage{amsmath, mathtools, amssymb, amsthm}
\usepackage{subfigure}
\usepackage[left=1in,right=1in,top=1in,bottom=1in]{geometry}
\usepackage{soul}
\usepackage{bm}
\usepackage{empheq}
\usepackage[normalem]{ulem}
\DeclareMathOperator*{\essinf}{Ess\ inf}

\usepackage[dvipsnames]{xcolor}

\numberwithin{equation}{section}

\renewcommand{\geq}{\geqslant}

\renewcommand{\leq}{\leqslant}

\newcommand{\tr}{\top}  

\makeatletter
\renewcommand*\env@matrix[1][\arraystretch]{%
  \edef\arraystretch{#1}%
  \hskip -\arraycolsep
  \let\@ifnextchar\new@ifnextchar
  \array{*\c@MaxMatrixCols c}}
\makeatother

\title[Multicomponent Euler equations]{A minimum entropy principle in the compressible multicomponent Euler equations}

\author[A. Gouasmi]{Ayoub Gouasmi}
\address{Ayoub Gouasmi\newline
\indent
Department of Aerospace Engineering, University of Michigan, Ann Arbor, MI}
\email{gouasmia@umich.edu}

\author[K. Duraisamy]{Karthik Duraisamy}
\address{Karthik Duraisamy\newline
\indent
Department of Aerospace Engineering, University of Michigan, Ann Arbor, MI}
\email{kdur@umich.edu}

\author[S. M. Murman]{Scott M. Murman}
\address{Scott M. Murman\newline
\indent
NASA Advanced Supercomputing Division, NASA Ames Research Center,  Moffett field, \newline \indent CA}
\email{Scott.M.Murman@nasa.gov}

\author[E. Tadmor]{Eitan Tadmor}
\address{Eitan Tadmor\newline
\indent
Department of Mathematics, Institute for Physical Sciences \& Technology (IPST)\newline
\indent and Center for Scientific Computation and Mathematical Modeling (CSCAMM)\newline
\indent
University of Maryland, College Park, MD}
\email{tadmor@cscamm.umd.edu}

\date{\today}

\begin{document}

\date{\today}



\thanks{\textbf{Acknowledgment.} AG and KD were funded by the AFOSR through grant number FA9550-16-1-030 (Tech. monitor: Fariba Fahroo). ET was supported in part by NSF grants DMS16-13911, RNMS11-07444 (KI-Net) and ONR grant N00014-1812465. 
 Ayoub Gouasmi is grateful to Jean-Luc Guermond, Bojan Popov and Ignacio Tomas for spirited discussions on their work \cite{Guermond_visc, Guermond_IDP, Guermond_IDP2} and for bringing the work of Delchini \textit{et al.} \cite{Delchini_1, Delchini_2} and Harten \textit{et al.} \cite{SuperHarten} to his attention.}

\begin{abstract}
 In this work, the space of admissible entropy functions for the compressible multicomponent Euler equations is explored, following up on [Harten, \textit{J. Comput. Phys.}, 49 (1), 1983, pp. 151-164]. This effort allows us to prove a minimum entropy principle on entropy solutions, whether smooth or discrete, in the same way it was originally demonstrated for the compressible Euler equations by [Tadmor, \textit{Appl. Numer. Math.}, 49 (3-5), 1986, pp. 211-219]. 
\end{abstract}

\maketitle

\section{Introduction}
Some hyperbolic systems of conservation laws,
\begin{equation}\label{eq:PDE_1D}
\partial_t \mathbf{u}  + \partial_x \mathbf{f} = 0,
\end{equation}
where $\mathbf{u}(x,t)$ and $\mathbf{f}(\mathbf{u}(x,t))$ are the state and flux vectors, respectively, admit a convex extension \cite{Friedrichs, Harten} in the sense that equation (\ref{eq:PDE_1D}) implies an additional conservation equation:
\begin{equation}\label{eq:Entropy}
    \partial_t U + \partial_x F = 0, 
\end{equation}
where $(U, \ F) = (U(\mathbf{u}), F(\mathbf{u})) \in \mathbb{R}^2$ is an entropy-entropy flux pair satisfying:
\begin{equation}\label{eq:entropy_cons}
    \frac{\partial U}{\partial \mathbf{u}}
    \frac{\partial \mathbf{f}}{\partial \mathbf{u}} = \frac{\partial F}{\partial \mathbf{u}}
\end{equation}
and $U$ strictly convex. We refer to $U$ as an \textit{entropy function}. Equation (\ref{eq:entropy_cons}) is a necessary and sufficient condition for (\ref{eq:PDE_1D}) to imply (\ref{eq:Entropy}). Additionally, Mock \cite{Mock} showed that the mapping $\mathbf{u} \rightarrow \mathbf{v}$ with the vector of entropy variables $\mathbf{v}$ defined as:
\begin{equation}
    \mathbf{v} := \bigg(\frac{\partial U}{\partial \mathbf{u}}\bigg)^{\tr},
\end{equation}
is one-to-one and turns (\ref{eq:PDE_1D}) into a symmetric hyperbolic system \cite{Friedrichs, Krushkov}. \\
\indent It is well known that when the flux $\mathbf{f}$ is nonlinear, discontinuous solutions to equation (\ref{eq:PDE_1D}) can develop from smooth initial conditions. Weak solutions must therefore be sought. Unfortunately, weak solutions are not uniquely defined and one needs additional conditions to distinguish physical solutions from non-physical ones. It is common practice to view physical solutions as those arising as vanishing viscosity limits, $\mathbf{u}(x,t) = \lim_{\epsilon \rightarrow 0} \mathbf{u}^{\epsilon}(x,t)$, of solutions $\mathbf{u}^{\epsilon}(x,t)$ to the regularized system:
\begin{equation}\label{eq:PDE_visc_1D}
\partial_t \mathbf{u}^{\epsilon}  + \partial_x \mathbf{f}(\mathbf{u}^{\epsilon}) = \epsilon \partial_x^2  \mathbf{u}^{\epsilon}, \ \epsilon > 0.
\end{equation}
Multiplying (\ref{eq:PDE_visc_1D}) on the left by $\mathbf{v}^{\tr}$ and using the convexity of $U$ one can show that $\mathbf{u}^{\epsilon}$ satisfies the inequality:
\begin{equation}
    \partial_t U(\mathbf{u}^{\epsilon}) + \partial_x F(\mathbf{u}^{\epsilon}) \leq \epsilon \partial_x^2  U(\mathbf{u}^{\epsilon}).
\end{equation}
In the limit $\epsilon \rightarrow 0$, this leads to the well-known entropy condition \cite{Lax, Krushkov}:
\begin{equation}\label{eq:entropy_cond}
    \partial_t U(\mathbf{u}) + \partial_x  F(\mathbf{u}) \leq 0,
\end{equation}
which is understood in the sense of distributions. Weak solutions to (\ref{eq:PDE_1D}) which satisfy the entropy condition (\ref{eq:entropy_cond}) \textit{for all entropies} are called \textit{entropy solutions}. \\
\indent For the compressible Euler equations governing the inviscid polytropic gas dynamics, Tadmor \cite{Tadmor2} showed that entropy solutions, whether smooth or discrete, satisfy a minimum entropy principle, namely that the spatial minimum of the specific entropy is an increasing function of time. \\
\indent In this work, we seek to extend this result to entropy solutions of the multicomponent compressible Euler equations. In section \ref{sec:GE}, we review the system at hand. In section \ref{sec:minS}, we recall the original proof and motivate the two families of entropy function we investigate in section \ref{sec:entropy_functions}. We end up showing a minimum entropy principle for the mixture's specific entropy. In section 5, we review numerical schemes which satisfy this property.  
\section{Governing equations}\label{sec:GE}
\hspace{0.4 cm} We consider the compressible multicomponent Euler equations \cite{Giovangigli} which consist of the conservation of species mass, momentum and total energy. In one dimension, that is equation (\ref{eq:PDE_1D}) with the state vector $\mathbf{u}$ and flux vector $\mathbf{f}$ defined by:
\begin{equation*}
    \mathbf{u} := \begin{bmatrix} \rho_1 & \hdots & \rho_N & \rho u & \rho e + \frac{1}{2} \rho u^2 \end{bmatrix}^{\tr}, \
    \mathbf{f} := \begin{bmatrix} \rho_1 u & \hdots & \rho_N u & \rho u^2 + p & (\rho e + \frac{1}{2} \rho u^2 + p)u \end{bmatrix}^{\tr},
\end{equation*}
where $\rho_k$ is the partial density of species $k$, $\rho := \sum_{k=1}^N \rho_k$ is the total density and $u$ is the fluid velocity. The pressure $p$ is given by the perfect gas law:
\begin{equation*}
p := \sum_{k=1}^N \rho_k  r_k T, \ r_k = \frac{R}{m_k},
\end{equation*}
where $m_k$ is the molar mass of species k and $R$ is the gas constant. The temperature $T$ is determined by the internal energy $\rho e$ which in this work is modeled following a thermally perfect gas assumption:
 \begin{equation*}
\rho e := \sum_{k=1}^N \rho_k e_k, \ e_k := e_{0k} + \int_{0}^{T} c_{vk}(\tau) d\tau.
\end{equation*}
For species k, $e_k$ is the specific internal energy of species k, $e_{0k}$ is a constant and $c_{vk} = c_{vk}(T) > 0$ is the constant volume specific heat.  Other quantities which will be used in this work are given by:
\begin{gather*}
h_k := e_k + r_k T, \
\rho c_v := \sum_{k=1}^N \rho_k c_{vk}, \ c_{pk} := c_{vk} + r_k, \ \gamma := \frac{c_p}{c_v}, \ Y_k := \frac{\rho_k}{\rho}.
\end{gather*}
$h_k$ is the specific enthalpy of species k, $c_v$ is the constant volume specific heat of the gas mixture, $\gamma$ is the specific heat ratio and $Y_k$ is the mass fraction of species $k$. The thermodynamic entropy of the mixture is given by:
\begin{gather*}
\rho s := \sum_{k=1}^N \rho_k s_k, \ 
s_k := \int_{0}^{T} \frac{c_{vk}(\tau)}{\tau} d\tau - r_k \ln(\rho_k)
\end{gather*}
Combining the transport equations for total density, species fractions and internal energy: 
\begin{equation}\label{eq:dRYE}
    D_t \rho = - \rho \partial_x u, \ D_t Y_k = 0, \ D_t e = -\frac{p}{\rho}\partial_x u,
\end{equation}
with the Gibbs relation:
\begin{equation}\label{eq:Gibbs0}
    T ds = de - \frac{p}{\rho^2} d\rho - \sum_{k=1}^N g_k dY_k,
\end{equation}
leads to a transport equation for the specific entropy $s$:
\begin{equation}\label{eq:entropy_transport}
    D_t s = 0.
\end{equation}
With total mass conservation, this leads to the conservation equation:
\begin{equation}\label{eq:S}
    \partial_t (\rho s) + \partial_x (\rho s u) = 0.
\end{equation}
For $\rho_k > 0, \ T > 0$, $(U, F) = (-\rho s, -\rho u s)$  is a valid entropy-entropy flux pair \cite{Chalot, Giovangigli}. The condition (\ref{eq:entropy_cons}) is met as a consequence of (\ref{eq:S}). The convexity of $U$ is established by looking at the entropy Hessian $\mathbf{G}$ given by:
\begin{equation*}
    \mathbf{G} := \frac{\partial^2 U}{\partial \mathbf{u}^2} = \frac{\partial \mathbf{v}}{\partial \mathbf{u}} = \frac{\partial \mathbf{v}}{\partial Z} \bigg( \frac{\partial \mathbf{u}}{\partial Z} \bigg)^{-1}.
\end{equation*}
The entropy variables $\mathbf{v}$ for the multicomponent system  can be easily derived using variable changes. Define the vector of primitive variables $Z = \begin{bmatrix} \rho_1 & \hdots & \rho_N & u & T \end{bmatrix}^{\tr}$. The chain rule gives:
\begin{equation*}
    \frac{\partial U}{\partial \mathbf{u}} = \frac{\partial U}{\partial Z} \bigg ( \frac{\partial \mathbf{u}}{\partial Z}\bigg)^{-1}.
\end{equation*}
The Gibbs identity (\ref{eq:Gibbs0}) can be written as:
\begin{equation}\label{eq:Gibbs}
    T dU = -d\rho e + \sum_{k=1}^N g_k d\rho_k,
\end{equation}
where $g_k =  h_k - T s_k$ is the Gibbs function of species k. From the definition of $\rho e$ we have:
\begin{equation}\label{eq:dE}
    d\rho e = \sum_{k=1}^N e_k d \rho_k + \rho c_v dT.
\end{equation}
Combining eqs. (\ref{eq:dE}) and (\ref{eq:Gibbs}), one obtains:
\begin{equation*}
    dU = \frac{1}{T} \bigg (\sum_{k=1}^N (g_k - e_k) d\rho_k - \rho c_v dT \bigg).
\end{equation*}
This gives:
\begin{equation}\label{eq:dSdZ}
    \frac{\partial U}{\partial Z} = \frac{1}{T} 
    \begin{bmatrix} (g_1 - e_1) & \hdots & (g_N - e_N) & 0 & -\rho c_v \end{bmatrix}.
\end{equation}
The Jacobian of the mapping $Z \rightarrow \mathbf{u}$ is given by:
\begin{equation}\label{eq:dudz}
    \frac{\partial \mathbf{u}}{\partial Z} = 
        \begin{bmatrix} 
            1 &        & 0 &    0   &    0   \\
              & \ddots &   & \vdots & \vdots \\
            0 &        & 1 &    0   &    0   \\
            u & \hdots & u &  \rho  &    0   \\
            e_{1} + k & \hdots & e_{N} + k & \rho u & \rho c_v 
        \end{bmatrix},
\end{equation}
where $k = \frac{1}{2}u^2$. The inverse of this matrix is given by:
\begin{equation}\label{eq:dudz_inv}
    \bigg(\frac{\partial \mathbf{u}}{\partial Z}\bigg)^{-1} = 
        \begin{bmatrix} 
            1 &        & 0 &    0   &    0   \\
              & \ddots &   & \vdots & \vdots \\
            0 &        & 1 &    0   &    0   \\
            -u \rho^{-1} & \hdots & -u \rho^{-1} &  \rho^{-1}  &    0   \\
            (k - e_{1})(\rho c_v)^{-1} & \hdots & (k - e_{N})(\rho c_v)^{-1} &  -u (\rho c_v)^{-1} & (\rho c_v)^{-1} 
        \end{bmatrix}.
\end{equation}
Combining eqs. (\ref{eq:dudz_inv}) and (\ref{eq:dSdZ}) yields the entropy variables \cite{Chalot, Giovangigli}:
\begin{equation}\label{eq:v}
    \mathbf{v} = \bigg(\frac{\partial U}{\partial \mathbf{u}}\bigg)^{\tr} = \frac{1}{T} 
    \begin{bmatrix} g_1 - k & \hdots & g_N - k & u & -1 \end{bmatrix}^{\tr}.
\end{equation}
We have:
\begin{equation}\label{eq:dvdZ}
    \frac{\partial \mathbf{v}}{\partial Z} = 
        \begin{bmatrix} 
            r_1/\rho_1  &        &     0      &    -u/T   &    (k - e_1) / T^2   \\
                        & \ddots &            & \vdots & \vdots \\
            0           &        & r_N/\rho_N &    -u/T   &    (k - e_N)/T^2   \\
            0           & \hdots &     0      &   1/T  & -u/T^2 \\
            0           & \hdots &     0      &    0   &  1/T^2 
        \end{bmatrix}.
\end{equation}
therefore the Hessian is given by:
\begin{equation}
    \mathbf{G} = \frac{1}{\rho c_v T^2} \begin{bmatrix}
                          &                      &                       & -u(k - (e_1 - c_v T)) & -(e_1 - k) \\
                          & \big(\zeta_{ij}\big) &                       &         \vdots        &   \vdots   \\
                          &                      &                       & -u(k - (e_N - c_v T)) & -(e_N - k) \\     
    -u(k - (e_1 - c_v T)) &       \hdots         & -u(k - (e_N - c_v T)) &      (u^2 + c_v T)    &     -u     \\
               -(e_1 - k) &       \hdots         &         -(e_N - k)    &            -u         &      1
    \end{bmatrix},
\end{equation}
with $ \zeta_{ij} = (\rho c_v T^2) \big(\delta_{ij} r_i/\rho_i + u^2 c_v T \big) + (e_i - k)(e_j - k)$ for $1 \leq i,j \leq N$. The positive definiteness of the Hessian matrix $\mathbf{G}$ is not immediately visible because it is dense. However the matrix $\mathbf{H}$ defined by the congruence relation:
\begin{equation}
    \mathbf{H} := \bigg( \frac{\partial \mathbf{u}}{\partial Z} \bigg)^{\tr} \mathbf{G} \bigg( \frac{\partial \mathbf{u}}{\partial Z} \bigg)  = \ \bigg( \frac{\partial \mathbf{u}}{\partial Z} \bigg)^{\tr} \frac{\partial \mathbf{v}}{\partial Z} = 
    \begin{bmatrix} 
            r_1/\rho_1  &        &     0      &    0   &    0   \\
                        & \ddots &            & \vdots & \vdots \\
            0           &        & r_N/\rho_N &    0   &    0   \\
            0           & \hdots &     0      & \rho/T &    0 \\
            0           & \hdots &     0      &    0   &  \rho c_v /T^2 
        \end{bmatrix},
\end{equation}
is positive definite, therefore $G$ is positive definite. This congruence relation, which was cleverly used in \cite{SuperHarten}, will be used as well in section \ref{sec:entropy_functions}. 
\section{The minimum entropy principle}\label{sec:minS} 
In this section, we review the proof of Tadmor \cite{Tadmor} for the compressible Euler equations then discuss how to apply it to the multicomponent system.

\subsection{Review}\label{sec:minS_proof}
Integrating the inequality (\ref{eq:entropy_cond}) over any domain $\Omega$ which induces no entropy influx across its boundaries gives:
\begin{equation}\label{eq:int_dSdt}
    \frac{d}{dt}\int_{\Omega} U (\mathbf{u}(x,t)) dx \leq 0 
\end{equation}
Integrating the above in time gives \cite{Lax}:
\begin{equation}
    \int_{\Omega} U (\mathbf{u}(x,t)) dx \leq  \int_{\Omega} U (\mathbf{u}(x,0)) dx
\end{equation}
Tadmor \cite{Tadmor0} showed that a sharper, more local version of the above inequality can be obtained:
\begin{equation}\label{eq:entropy_sharp}
    \int_{|x| \leq R} U (\mathbf{u}(x,t)) dx \leq  \int_{|x| \leq R + t \cdot q_{max}} U (\mathbf{u}(x,0)) dx,
\end{equation}
where $q_{max}$ is the maximum velocity in the domain at $t = 0$. For the Euler equations, Harten \cite{Harten} sought pairs of the form $(U^h, F^h) = (- \rho h(s), - \rho u h(s))$ where $s = \ln (p) - \gamma \ln(\rho)$ is the dimensionless specific entropy (divided by the $c_v$, we will use the letter $f$ instead of $h$ in section \ref{sec:entropy_functions}) and $h$ is a smooth function of $S$. Harten showed that the pair $(U^{h}, F^{h})$ is admissible if and only if $h$ satisfies:
\begin{equation}\label{eq:Harten_cond}
    h^{'} - \gamma \ h^{''} > 0, \ h^{'} > 0.
\end{equation}
For any such function $h$, the inequality (\ref{eq:entropy_sharp}) with $U = U^{h}$ gives:
\begin{equation}\label{eq:entropy_sharp2}
    \int_{|x| \leq R} \rho (x, t) \cdot h(s(x,t)) \ dx \geq  \int_{|x| \leq R + t \cdot q_{max}} \rho (x, 0) \cdot h(s(x,0)) \ dx.
\end{equation}
Tadmor makes a special choice $h_0$ for the function $h$:
\begin{equation*}
    h_0(s) = \min[ s - s_0, \ 0], \ s_0 = \essinf_{|x| \leq R + t \cdot q_{max}} s(x, 0).
\end{equation*}
$s_0$ is the essential infimum of the specific entropy in the domain $\Omega = \{ x : |x| < R + t \cdot q_{max} \}$. From inequality (\ref{eq:entropy_sharp2}), we get:
\begin{equation}\label{eq:entropy_sharp3}
    \int_{|x| \leq R} \rho (x, t) \cdot \min[s(x,t) - s_0, \ 0] \ dx \geq  \int_{|x| \leq R + t \cdot q_{max}} \rho (x, 0) \cdot \min[s(x,0) - s_0, \ 0] \ dx.
\end{equation}
The right-hand side drops by definition of $s_0$, so equation (\ref{eq:entropy_sharp3}) simplifies to:
\begin{equation}\label{eq:entropy_sharp4}
    \int_{|x| \leq R} \rho (x, t) \cdot \min[s(x,t) - s_0, \ 0] \ dx \geq  0.
\end{equation}
The integrand on the left-hand side is negative, therefore inequality (\ref{eq:entropy_sharp4}) imposes for $|x| \leq R$:
\begin{equation}\label{eq:min_entropy_OP}
    \min[s(x,t) - s_0, \ 0] = 0 \Leftrightarrow  s(x,t) \geq  \essinf_{|x| \leq R + t \cdot q_{max}} s(x, 0).
\end{equation}
This is the minimum entropy principle satisfied by \textit{entropy solutions to the compressible Euler equations}. A similar result holds for discrete solutions $\mathbf{u}_{i}^{n}$ (the subscript $i$ and the superscript $n$ refer to the cell index and time instant, respectively) which satisfy the fully-discrete entropy inequality:
\begin{equation}\label{eq:discrete_ES}
    \sum_{i} U(\mathbf{u}_i^{n+1}) \leq \sum_{i} U(\mathbf{u}_i^{n}),
\end{equation}
for all entropies $U$. Taking $U = - \rho h_0(s)$ with $s_0$ defined as the minimum specific entropy at time instant $n$ leads to:
\begin{equation*}
    \sum_{i} \rho(\mathbf{u}_i^{n+1}) \cdot \min[s(\mathbf{u}_{i}^{n+1}) - s_0, \ 0] \geq 0.
\end{equation*}
If $\rho(\mathbf{u}_i^{n+1}) > 0$, this imposes in every cell:
\begin{equation}\label{eq:min_entropy}
    \min[s(\mathbf{u}_{i}^{n+1}) - s_0, \ 0] = 0 \ \Leftrightarrow s(\mathbf{u}_{i}^{n+1}) \geq  \min_{i} s(\mathbf{u}_i^{n}).
\end{equation}
\indent At first glance, injecting $U = -\rho h_0(s)$ in inequalities (\ref{eq:entropy_sharp}) and (\ref{eq:discrete_ES}) should not be allowed because $h_0$ is not smooth function of $s$. What makes this step valid nonetheless is the fact that $h_0$ \textit{can be written as the limit of a sequence of smooth functions which satisfy Harten's conditions}. Without loss of generality, let's assume $s_0 = 0$ and consider the convolution defined as:
\begin{equation*}
    h(s) = \int_{-\infty}^{+\infty} h_0(s - \overline{s}) \phi(\overline{s}) d\overline{s}.
\end{equation*}
where $\phi$ is a smooth function satisfying:
\begin{equation*}
    \int_{-\infty}^{+\infty} \phi(\overline{s})d\overline{s} = 1, \ \phi(\overline{s}) > 0.
\end{equation*}
$\phi$ should also be such that the convolution is well-defined everywhere. $\phi(\overline{s}) = \exp(-\overline{s}^2)/\sqrt{\pi}$ is a valid choice. By definition of $h_0$, we have:
\begin{equation*}
    h(s) = \int_{s}^{+\infty} (s - \overline{s}) \phi(\overline{s}) d\overline{s} = s \int_{s}^{+\infty} \phi(\overline{s})d\overline{s} - \int_{s}^{+\infty} \overline{s}\phi(\overline{s}) d\overline{s}.
\end{equation*}
$h$ is smooth and satisfies Harten's conditions because:
\begin{equation*}
    h^{'}(s) = \int_{s}^{+\infty} \phi(\overline{s})d\overline{s} > 0, \ h^{''}(s) = - \phi(s) < 0.
\end{equation*}
$\forall \varepsilon > 0$, the function $h_{\varepsilon}$ defined by:
\begin{equation}\label{eq:convo_H}
    h_{\varepsilon}(s) = \int_{-\infty}^{+\infty} h_0(s - \overline{s}) \phi_{\varepsilon}(\overline{s}) d\overline{s}, \ \phi_{\varepsilon}(\overline{s}) = \frac{1}{\varepsilon} \phi\bigg(\frac{\overline{s}}{\varepsilon}\bigg),
\end{equation}
is smooth and satisifies Harten's conditions as well. What is more, $\phi_{\varepsilon}$ converges, in the sense of distributions, to the Dirac delta function when $\varepsilon \rightarrow 0$ (classic result). Therefore, inequality (\ref{eq:entropy_sharp3}) is obtained  $h_0 = \lim_{\varepsilon \rightarrow} h_{\varepsilon}$. \\
\indent The main takeaway of this review is that \textit{not all entropy inequalities need to be satisfied for a minimum entropy principle to hold in the compressible Euler equations}. Those involving the "convolution entropies" $U = - \rho h_{\varepsilon}(s), \forall \varepsilon > 0$ defined by equation (\ref{eq:convo_H}) are enough to conclude. \\ \\
\indent \ul{Remark 1}: This proof and Harten's characterization (\ref{eq:Harten_cond}) are both independent of the number of spatial dimensions \cite{Harten, Tadmor2}. Throughout this manuscript, we are working in one dimension for the sake of simplicity only. \\
\indent \ul{Remark 2}: Kroner \textit{et al.} \cite{Kroner} use a different approach to demonstrate that bounded entropy solutions to the quasi-1D Euler equations with discontinuous cross-section satisfy a minimum entropy principle. The inequality (\ref{eq:entropy_sharp2}) is used with $h(s) = -(s_0 - s)^p, \ p > 1, s_0 > s$ ($s_0$ denotes an upper bound in this context), raised to the power $1/p$ and passed to the limit $p \rightarrow \infty$. \\
\indent \ul{Remark 3}: A minimum entropy principle for smooth solutions to well-designed regularizations of the Euler equations was proved by Guermond and Popov \cite{Guermond_visc} (see also Delchini \textit{et al.} \cite{Delchini_1, Delchini_2} for other systems). In this work, we are interested in the minimum entropy principle as a property of entropy solutions, \textit{whether smooth or discrete}, to the multicomponent compressible Euler equations. 

\subsection{Elements of proof for the multicomponent compressible Euler equations}
\hspace*{0.4 cm} We need to formulate \textit{what a minimum entropy principle would be} in the multicomponent case. The first option is a minimum entropy principle \textit{involving the specific entropy of each species}:
\begin{equation*}
    s_k(x,t) \geq  s_{0k} = \essinf_{|x| \leq R + t \cdot q_{max}} s_k(x, 0), \ 1 \leq k \leq N.
\end{equation*}
Working Tadmor's proof backwards, this is obtained if we can show that entropy solutions satisfy the inequality:
\begin{equation}\label{eq:entropy_MC_sharp2}
    \int_{|x| \leq R} \sum_{k=1}^N \rho_k (x, t) \cdot f_k(s_k(x,t)) \ dx \geq  \int_{|x| \leq R + t \cdot q_{max}}  \sum_{k=1}^N \rho_k (x, 0) \cdot f_k(s_k(x,0)) \ dx,
\end{equation}
and that $f_k$ can be taken as $f_{0k}(s_k) = \min[s_k - s_{0k}, \ 0]$. This leads us to examine entropy pairs $(U_I^{f}, F_I^{f})$ of the form:
\begin{equation}\label{eq:entropy_1}
    (U_I^{f}, F_I^{f}) = \bigg(- \sum_{k=1}^N \rho_k f_k, \ - \sum_{k=1}^N \rho_k u f_k \bigg), \ f_k = f_k(s_k),
\end{equation}
and attempt to show that those with $f_k$ defined as the convolution  (\ref{eq:convo_H}) are valid entropy pairs. The second option is a minimum entropy principle \textit{involving the specific entropy of the gas mixture}:
\begin{equation*}
    s(x,t) \geq  s_0 = \essinf_{|x| \leq R + t \cdot q_{max}} s(x, 0).
\end{equation*}
In the same vein, this is obtained if we can show that entropy solutions satisfy the inequality:
\begin{equation}\label{eq:entropy_MC_sharp3}
    \int_{|x| \leq R} \rho (x, t) \cdot f(s(x,t)) \ dx \geq  \int_{|x| \leq R + t \cdot q_{max}} \rho (x, 0) \cdot f(s(x,0)) \ dx,
\end{equation}
and that $f$ can be taken as $f_0(s) = \min[s - s_{0}, \ 0]$. This leads us to examine entropy pairs $(U_{II}^f, F_{II}^f)$ of the form:
\begin{equation}\label{eq:entropy_2}
    (U_{II}^f, F_{II}^f) = (-\rho f(s), -\rho u f(s)), 
\end{equation}
and attempt show that those with $f$ defined as the convolution (\ref{eq:convo_H}) are valid entropy pairs. \\
\indent These two families are investigated in the next section. The admissibility conditions will take the form of constraints of the first and second derivatives of $f_k$ (first case) and $f$ (second case). If the first and second derivatives are allowed to be strictly positive and negative, respectively, then the convolution (\ref{eq:convo_H}) qualifies and a minimum entropy principle follows.

\section{Entropy functions in the multicomponent case}\label{sec:entropy_functions}
For each candidate family of entropy functions, we must check for conservation and convexity with respect to the conservative variables. For a candidate entropy $U^f$, convexity is equivalent to the positive definiteness of its Hessian matrix $\mathbf{G}$:
\begin{equation*}
    \mathbf{G} = \frac{\partial^2 U^{f}}{\partial \mathbf{u}^2} = \frac{\partial \mathbf{v}^f}{\partial \mathbf{u}}, \ \mathbf{v}^f = \bigg( \frac{\partial U^{f}}{\partial \mathbf{u}}\bigg)^{\tr}.
\end{equation*}
$\mathbf{v}^f$ is the vector of entropy variables associated with the candidate entropy.
\subsection{Candidate I}

\subsubsection*{Conservation}
Equation (\ref{eq:Entropy}) with $(U, F) = (U_I^{f}, F_I^{f})$ holds if and only if $\sum_{k=1}^N Y_k f_k$ satisfies a transport equation. We have:
\begin{align*}
    d\bigg(\sum_{k=1}^N Y_k f_k\bigg) =& \ \sum_{k=1}^N Y_k df_k + \sum_{k=1}^N f_k dY_k \\
                                =& \ \sum_{k=1}^N Y_k f_k^{'} ds_k + \sum_{k=1}^N f_k dY_k \\
                                =& \ \sum_{k=1}^N Y_k f_k^{'} \bigg(\frac{c_{vk}}{T}dT - \frac{r_k}{\rho_k} d\rho_k \bigg) + \sum_{k=1}^N f_k dY_k \\
                                =& \ \bigg(\sum_{k=1}^N Y_k f_k^{'} c_{vk}\bigg) \frac{dT}{T} - \frac{1}{\rho}\sum_{k=1}^N f_k^{'} r_k d\rho_k + \sum_{k=1}^N f_k dY_k \\
                                =& \ \bigg(\sum_{k=1}^N Y_k f_k^{'} c_{vk}\bigg) \frac{dT}{T} - \bigg(\sum_{k=1}^N f_k^{'} Y_k r_k\bigg) \frac{d\rho}{\rho}  + \sum_{k=1}^N (f_k - r_k f_k^{'}) dY_k.
\end{align*}
From the differential relation:
\begin{equation*}
    de = \sum_{k=1}^N dY_k e_k + \sum_{k=1}^N Y_k c_{vk} dT = \sum_{k=1}^N dY_k e_k + c_v dT,
\end{equation*}
we obtain the following equation for temperature:
\begin{equation}\label{eq:dT}
    D_t T = -\frac{p}{\rho c_v} \partial_x u = \frac{p}{\rho^2 c_v} D_t \rho.
\end{equation}
Using equations (\ref{eq:dRYE}) and (\ref{eq:dT}), we can show that $U_{I}^f$ is conserved if and only if:
\begin{equation}
    \frac{1}{T}\bigg(\sum_{k=1}^N Y_k f_k^{'} c_{vk}\bigg) D_t T - \frac{1}{\rho}\bigg(\sum_{k=1}^N f_k^{'} Y_k r_k\bigg) D_t \rho = 0 \ \Leftrightarrow \frac{p}{\rho T} \bigg(\frac{\sum_{k=1}^N Y_k f_k^{'} c_{vk}}{\sum_{k=1}^N Y_k c_{vk}}\bigg)  - \bigg(\sum_{k=1}^N f_k^{'} Y_k r_k\bigg) = 0
\end{equation}
Using the ideal gas law, this condition rewrites:
\begin{equation}\label{eq:conservation}
    \frac{\sum_{k=1}^N \rho_k c_{vk} f_k^{'}}{\sum_{k=1}^N \rho_k c_{vk}} = \frac{\sum_{k=1}^N \rho_k r_{k} f_k^{'}}{\sum_{k=1}^N \rho_k r_{k}}.
\end{equation}
\subsubsection*{Convexity}
We have:
\begin{equation*}
    \frac{\partial s_k}{\partial \rho_k} = -\frac{r_k}{\rho_k}, \ \ \frac{\partial s_k}{\partial T} = \frac{c_{vk}}{T}, \ \ \frac{\partial f_k}{\partial \rho_k} = -\frac{r_k}{\rho_k} f_k^{'}, \ \ \frac{\partial f_k}{\partial T} = \frac{c_{vk}}{T} f_k^{'}.
\end{equation*}
Therefore
\begin{equation*}
    \frac{\partial U_I^{f}}{\partial Z} = \begin{bmatrix}
            - f_1 + r_1 f_1^{'} & \dots & - f_N + r_N f_N^{'} & 0 & - \frac{1}{T}\big( \sum_{k=1}^N \rho_k c_{vk} f_k^{'} \big)
        \end{bmatrix},
\end{equation*}
and the entropy variables (chain rule) are given by:
\begin{equation*}
    \mathbf{v}_{I}^f = \begin{bmatrix}
            - f_1 + r_1 f_1^{'} - \beta \frac{k - e_1}{T} & \dots & - f_N + r_N f_N^{'} - \beta \frac{k - e_N}{T} & \beta \frac{u}{T} & - \beta \frac{1}{T} 
        \end{bmatrix}^{\tr}, \ \ \beta = \frac{ \sum_{k=1}^N \rho_k c_{vk} f_k^{'} }{\sum_{k=1}^N \rho_k c_{vk}}.
\end{equation*}
For simplicity, let's assume calorically perfect gases ($c_{vk}$ and $c_{pk}$ constants) and drop the standard formation constants. To proceed with the Hessian calculation we need the following:
\begin{equation*}
    \frac{\partial \beta}{\partial \rho_k} = \frac{c_{vk}}{\rho c_v} (f_k^{'} - r_k f_k^{''} - \beta), \ \ \frac{\partial \beta}{\partial T} = \frac{\eta}{T}, \ \eta = \frac{ \sum_{k=1}^N \rho_k c_{vk}^2 f_k^{''}}{\sum_{k=1}^N \rho_k c_{vk}}.
\end{equation*}
Denote $\xi_k = f_k^{'} - r_k f_k^{''}$ and $\mathbf{v}_I^f = [v_{1,1}^f \ \dots \ v_{1, N}^f  \ v_2^f \ v_3^f]^{\tr}$. The gradients of the last component are given by:
\begin{equation*}
    \frac{\partial v_3^f}{\partial \rho_k} = -\frac{1}{T}\frac{c_{vk}}{\rho c_v} (\xi_k - \beta), \
    \frac{\partial v_3^f}{\partial u} = 0, \ 
    \frac{\partial v_3^f}{\partial T} = \frac{\beta - \eta}{T^2}.
\end{equation*}
The gradients of the before-last component are given by:
\begin{equation*}
    \frac{\partial v_2^f}{\partial \rho_k} = \frac{u}{T} \frac{c_{vk}}{\rho c_v} (\xi_k - \beta), \
    \frac{\partial v_2^f}{\partial u} = \frac{\beta}{T}, \ 
    \frac{\partial v_2^f}{\partial T} = u \frac{\eta - \beta}{T^2}.
\end{equation*}
The gradient of the $l$-th component is given by:
\begin{equation*}
    \frac{\partial v_{1,l}^f}{\partial \rho_k} =  \delta_{kl} \frac{r_k}{\rho_k}\xi_k  - (\frac{k}{T} - c_{vl}) \frac{c_{vk}}{\rho c_v} (\xi_k - \beta), \
    \frac{\partial v_{1,l}^f}{\partial u} = - u \frac{\beta}{T}, \ 
    \frac{\partial v_{1,l}^f}{\partial T} = -\frac{c_{vl}}{T} \xi_l  + \frac{(\beta - \eta)k}{T^2} + c_{vl}\frac{\eta}{T}.
\end{equation*}
For two species, we have:
\begin{equation}\label{eq:dvfvz}
    \frac{\partial \mathbf{v}_{I}^f}{\partial Z} = \begin{bmatrix}
            \frac{r_1}{\rho_1}\xi_1 - (\frac{k}{T} - c_{v1})\frac{c_{v1}}{\rho c_v}(\xi_1 - \beta) & - (\frac{k}{T} - c_{v1}) \frac{c_{v2}}{\rho c_v} (\xi_2 - \beta) & -u \frac{\beta}{T} & -\frac{c_{v1}}{T} \xi_1  + \frac{(\beta - \eta)k}{T^2} + c_{v1}\frac{\eta}{T} \\
            - (\frac{k}{T} - c_{v2}) \frac{c_{v1}}{\rho c_v} (\xi_1 - \beta) & \frac{r_2}{\rho_2}\xi_2 - (\frac{k}{T} - c_{v2})\frac{c_{v2}}{\rho c_v}(\xi_2 - \beta) & -u \frac{\beta}{T} & -\frac{c_{v2}}{T} \xi_2  + \frac{(\beta - \eta)k}{T^2} + c_{v2}\frac{\eta}{T} \\
            \frac{u}{T} \frac{c_{v1}}{\rho c_v} (\xi_1 - \beta) & \frac{u}{T} \frac{c_{v2}}{\rho c_v} (\xi_2 - \beta) & \frac{\beta}{T} & u \frac{\eta - \beta}{T^2} \\
            -\frac{1}{T}\frac{c_{v1}}{\rho c_v} (\xi_1 - \beta) & -\frac{1}{T}\frac{c_{v2}}{\rho c_v} (\xi_2 - \beta) & 0 & \frac{\beta - \eta}{T^2}
        \end{bmatrix}.
\end{equation}
If $f(s) = s$ then $\beta = 1, \eta = 0$ and $\xi_k = 1$ and equation (\ref{eq:dvfvz}) does simplify to equation (\ref{eq:dvdZ}). The chain rule gives for the Hessian $\mathbf{G}_I$:
\begin{equation*}
    \mathbf{G}_I = \frac{\partial \mathbf{v}_I^f}{\partial Z} \bigg( \frac{\partial \mathbf{u}}{\partial Z} \bigg)^{-1}.
\end{equation*}
$\mathbf{G}_I$ is dense. We establish conditions on $f_k$ so that $G$ is positive definite by looking at the congruent matrix:
\begin{equation*}
    \mathbf{H}_I = \bigg( \frac{\partial \mathbf{u}}{\partial Z} \bigg)^{\tr} \mathbf{G}_I \bigg( \frac{\partial \mathbf{u}}{\partial Z} \bigg)  = \ \bigg( \frac{\partial \mathbf{u}}{\partial Z} \bigg)^{\tr} \frac{\partial \mathbf{v}_I^f}{\partial Z}.
\end{equation*}
$\mathbf{H}_I$ is given by:
\[
\begin{split}
    \mathbf{H}_I & = \ \begin{bmatrix} 
            1 & 0 &    u   &    c_{v1} T + k\\
            0 & 1 &    u   &    c_{vN} T + k  \\
            0 &  0 &  \rho  &    \rho u   \\
            0 &  0 &    0   & \rho c_v 
        \end{bmatrix}\\
        & \qquad  \times  
        \begin{bmatrix}
            \frac{r_1}{\rho_1}\xi_1 - (\frac{k}{T} - c_{v1})\frac{c_{v1}}{\rho c_v}(\xi_1 - \beta) & - (\frac{k}{T} - c_{v1}) \frac{c_{v2}}{\rho c_v} (\xi_2 - \beta) & -u \frac{\beta}{T} & -\frac{c_{v1}}{T} \xi_1  + \frac{(\beta - \eta)k}{T^2} + c_{v1}\frac{\eta}{T} \\
            - (\frac{k}{T} - c_{v2}) \frac{c_{v1}}{\rho c_v} (\xi_1 - \beta) & \frac{r_2}{\rho_2}\xi_2 - (\frac{k}{T} - c_{v2})\frac{c_{v2}}{\rho c_v}(\xi_2 - \beta) & -u \frac{\beta}{T} & -\frac{c_{v2}}{T} \xi_2  + \frac{(\beta - \eta)k}{T^2} + c_{v2}\frac{\eta}{T} \\
            \frac{u}{T} \frac{c_{v1}}{\rho c_v} (\xi_1 - \beta) & \frac{u}{T} \frac{c_{v2}}{\rho c_v} (\xi_2 - \beta) & \frac{\beta}{T} & u \frac{\eta - \beta}{T^2} \\
            -\frac{1}{T}\frac{c_{v1}}{\rho c_v} (\xi_1 - \beta) & -\frac{1}{T}\frac{c_{v2}}{\rho c_v} (\xi_2 - \beta) & 0 & \frac{\beta - \eta}{T^2}
        \end{bmatrix}\\
               & = \  \begin{bmatrix}
               \frac{r_1}{\rho_1}\xi_1 & 0 & 0 & -\frac{c_{v1}}{T}(\xi_1 - \beta) \\
               0 & \frac{r_2}{\rho_2}\xi_2 & 0 & -\frac{c_{v2}}{T}(\xi_2 - \beta) \\
               0 & 0 & \frac{\rho \beta}{T} & 0 \\
               -\frac{c_{v1}}{T}(\xi_1 - \beta) & -\frac{c_{v2}}{T}(\xi_2 - \beta) & 0 & \rho c_v \frac{\beta - \eta}{T^2} 
               \end{bmatrix}
\end{split}
\]
$\mathbf{H}_I$ is positive definite if and only if the determinants of the major blocks of $\mathbf{H}_I$ are all positive (from Harten \cite{Harten}). For the first three major blocks, this is equivalent to the requirement that $\xi_1 > 0$, $\xi_2 > 0$ and $\beta > 0$ are positive. Last:
\begin{align*}
    det(\mathbf{H}_I) =& \ \frac{\rho \beta}{T^3} r_1 r_2 \bigg( \rho c_v (\beta - \eta)\frac{\xi_1 \xi_2}{\rho_1 \rho_2 }  - \frac{c_{v1}}{\gamma_1-1}(\xi_1 - \beta)^2  \frac{\xi_2}{\rho_2} - \frac{c_{v2}}{\gamma_2-1}(\xi_2 - \beta)^2  \frac{\xi_1}{\rho_1}  \bigg) \\
                    =& \ \frac{\rho \beta r_1 r_2 \xi_1 \xi_2}{\rho_1 \rho_2 T^3}  \bigg( \rho c_v (\beta - \eta)  - \frac{\rho_1 c_{v1}}{\gamma_1-1}\frac{(\xi_1 - \beta)^2}{\xi_1}   - \frac{\rho_2 c_{v2}}{\gamma_2-1}\frac{(\xi_2 - \beta)^2}{\xi_2}  \bigg) \\
                     =& \ \frac{\rho \beta r_1 r_2 \xi_1 \xi_2}{\rho_1 \rho_2 T^3}  \bigg( \rho_1 c_{v1} \bigg((\beta - \eta)  - \frac{1}{\gamma_1-1}\frac{(\xi_1 - \beta)^2}{\xi_1}\bigg) + \rho_2 c_{v2} \bigg((\beta - \eta)  - \frac{1}{\gamma_2-1}\frac{(\xi_2 - \beta)^2}{\xi_2}\bigg)  \bigg) \\
                     =& \ \frac{\rho \beta r_1 r_2 \xi_1 \xi_2}{\rho_1 \rho_2 T^3}  \bigg( \frac{\rho_1 c_{v1}}{\xi_1(\gamma_1-1)} \bigg((\beta - \eta)\xi_1(\gamma_1-1)  - (\xi_1 - \beta)^2\bigg) \\
                     & \quad + \frac{\rho_2 c_{v2}}{\xi_2(\gamma_2-1)} \bigg((\beta - \eta)\xi_2(\gamma_2-1)  - (\xi_2 - \beta)^2\bigg)  \bigg) \\ 
                     =& \ \frac{\rho \beta r_1 r_2 \xi_1 \xi_2}{\rho_1 \rho_2 T^3}  \bigg( \frac{\rho_1 c_{v1}}{\xi_1(\gamma_1-1)} \Delta_1 + \frac{\rho_2 c_{v2}}{\xi_2(\gamma_2-1)} \Delta_2  \bigg),
\end{align*}
where $\Delta_k = (\beta - \eta)\xi_k(\gamma_k-1)  - (\xi_k - \beta)^2$. For an arbitrary number of species:
\begin{equation}\label{eq:H1}
    \mathbf{H}_I = \begin{bmatrix}
                  \frac{r_1}{\rho_1}\xi_1       &        &                                  &          0            & -\frac{c_{v1}}{T}(\xi_1 - \beta) \\
                                                & \ddots &                                  &       \vdots          &             \vdots               \\
                                                &        &     \frac{r_N}{\rho_N}\xi_N      &          0            & -\frac{c_{vN}}{T}(\xi_N - \beta) \\
                          0                     & \hdots &                0                 & \frac{\rho \beta}{T}  &                0                 \\
               -\frac{c_{v1}}{T}(\xi_1 - \beta) & \hdots & -\frac{c_{vN}}{T}(\xi_N - \beta) &          0            & \rho c_v \frac{\beta - \eta}{T^2} 
               \end{bmatrix},
\end{equation}
and one can easily show that:
\begin{equation}
    det(\mathbf{H}_I) = \frac{\rho \beta}{T^3} \bigg(\prod_{k=1}^N \frac{r_k \xi_k}{\rho_k} \bigg) \bigg( \sum_{k=1}^N \frac{\rho_k c_{vk}}{\xi_k(\gamma_k - 1)} \Delta_k \bigg).
\end{equation}
Overall, $U^f$ is an admissible entropy for the multicomponent Euler equations if and only if: 
\begin{equation}\label{eq:Harten_MC}
    \frac{\sum_{k=1}^N \rho_k c_{vk} f_k^{'}}{\sum_{k=1}^N \rho_k c_{vk}} = \frac{\sum_{k=1}^N \rho_k r_{k} f_k^{'}}{\sum_{k=1}^N \rho_k r_{k}}, \ \  \xi_k > 0, \ \ \beta > 0, \ \ \sum_{k=1}^N \frac{\rho_k c_{vk}}{\xi_k(\gamma_k - 1)} \Delta_k > 0.
\end{equation}
While the sufficient conditions $f_k^{'} > 0, \ f_k^{''} < 0$ for a minimum entropy principle are compatible with $\xi_k > 0$ and $\beta > 0$, it is not clear whether they are compatible with the last inequality of (\ref{eq:Harten_MC}) ($\Delta_k$ being the difference of two positive terms). Additionally, the equality constraint (\ref{eq:conservation}) which came from the requirement of conservation does not seem to offer any option other than $f_k^{'}$ constant. Note that if $f_k^{'} > 0, \ f_k^{''} < 0$ were to violate any of the conditions derived here, it would only mean that we cannot prove a minimum entropy principle with the approach exposed in section \ref{sec:minS_proof}. Disproving a minimum entropy principle would require a counterexample. \\
\indent For the compressible Euler equations, $\mathbf{H}_I$ simplifies to:
\begin{equation*}
    \mathbf{H}_I = \begin{bmatrix}
               \frac{r}{\rho}\xi & 0 & -\frac{c_{v}}{T}(\xi - \beta) \\
               0 & \frac{\rho \beta}{T} & 0 \\
               -\frac{c_{v}}{T}(\xi - \beta) & 0 & \rho c_v \frac{\beta - \eta}{T^2} 
               \end{bmatrix}, \ \xi = f^{'} - r f^{''}, \ \beta = f^{'}, \ \eta = c_v f^{''}.
\end{equation*}
The determinants of the three major blocks are:
\begin{equation*}
    det(H_{11}) = \frac{r}{\rho} \xi, \ det(H_{22}) = \frac{\rho }{T}\beta, \ \det(\mathbf{H}_I) =\frac{\rho r c_v \beta}{T^3(\gamma-1)} \bigg( (\beta - \eta)\xi (\gamma - 1) -  (\xi - \beta)^2\bigg).
\end{equation*}
Using $(\gamma -1 )(\beta - \eta) = (\gamma - 1)f^{'} - r f^{''}$ and $\xi - \beta = - rf^{''}$, the determinant simplifies to:
\begin{equation*}
    det(\mathbf{H}_I) = \frac{\rho r c_v \beta^2}{T^3} \bigg(  f^{'} - c_p f^{''} \bigg)  
\end{equation*}
The necessary conditions for $\mathbf{H}_I$ to be positive definite are then:
\begin{equation}\label{eq:Harten_0}
    f^{'} - r f^{''} > 0, \ f^{'} > 0, \ f^{'} - c_p f^{''} > 0.
\end{equation}
Since $f^{'} > 0$, the first and third inequality of (\ref{eq:Harten_0}) can be rewritten as:
\begin{equation*}
    \frac{f^{''}}{f^{'}} < \frac{1}{r}, \ \frac{f^{''}}{f^{'}} < \frac{1}{c_p}.
\end{equation*}
Since $c_p > r$, the first inequality is implied by the second. Therefore, the necessary conditions (\ref{eq:Harten_0}) simplify to:
\begin{equation}
    f^{'} > 0, \ f^{'} - c_{p} f^{''} > 0.
\end{equation}
These are the well-known conditions (\ref{eq:Harten_cond}) for the Euler equations (note that the function $f$ in this section and the function $h$ in section \ref{sec:minS_proof} are related by $f(s) = h(s/c_v)$). The conditions (\ref{eq:Harten_MC}) are therefore consistent with Harten's in the Euler case.
\subsection{Candidate II}
\subsubsection*{Conservation}
Multiplying the transport equation for the specific entropy (\ref{eq:entropy_transport}) with $f^{'}$ leads to a transport equation for $f(s)$. Conservation of $U_{II}^{f}$ with the entropy flux $F_{II}^{f}$ then follows from the total mass conservation equation.
\subsubsection*{Convexity}
We have:
\begin{equation*}
    \frac{\partial Y_j}{\partial \rho_k} = \frac{\delta_{jk}}{\rho} - \frac{\rho_j}{\rho^2}, \ \frac{\partial s}{\partial \rho_k} = \frac{1}{\rho}(s_k - r_k - s), \ \frac{\partial s}{\partial T} = \frac{c_v}{T}.
\end{equation*}
This gives:
\begin{equation}
    \frac{\partial U_{II}^f}{\partial Z} = \begin{bmatrix}
        f^{'} (- s_1 + r_1 + s) - f & \dots & f^{'} (- s_N + r_N + s) - f & 0 & -\frac{\rho c_v}{T} f^{'}
    \end{bmatrix},
\end{equation}
and the entropy variables:
\begin{equation}
    \mathbf{v}_{II}^f = \begin{bmatrix}
        f^{'} \frac{g_1 - k}{T}  + f^{'}s - f & \dots &  f^{'} \frac{g_N - k}{T}  + f^{'}s - f & f^{'} \frac{u}{T} & -f^{'}\frac{1}{T} 
    \end{bmatrix}^{\tr}  =  f^{'} \mathbf{v} + (f^{'} s - f ) \begin{bmatrix}
        1 & \cdots & 1 & 0 & 0
    \end{bmatrix}^{\tr}.
\end{equation}
Again, the conditions for convexity are established by looking at the congruent matrix $\mathbf{H}_{II}$ defined by:
\begin{equation*}
     \mathbf{H}_{II} = \bigg( \frac{\partial \mathbf{u}}{\partial Z} \bigg)^{\tr} \mathbf{G}_{II} \bigg( \frac{\partial \mathbf{u}}{\partial Z} \bigg)  = \ \bigg( \frac{\partial \mathbf{u}}{\partial Z} \bigg)^{\tr} \frac{\partial \mathbf{v}_{II}^f}{\partial Z}.
\end{equation*}
We have:
\begin{equation*}
        \frac{\partial \mathbf{v}_{II}^f}{\partial Z} = f^{'} \frac{\partial \mathbf{v}}{\partial Z} + \frac{f^{''}}{\rho} \begin{bmatrix} (g_1 - k)/T + s \\ \vdots \\ (g_N - k)/T + s \\ u/T \\ -1/T \end{bmatrix} 
            \begin{bmatrix} s_1 - r_1 - s & \hdots & s_N - r_N - s & 0 & \frac{\rho c_v}{T} \end{bmatrix} 
\end{equation*}
and 
\begin{equation*}
    \bigg(\frac{\partial \mathbf{u}}{\partial Z}\bigg)^{\tr}
    \begin{bmatrix} (g_1 - k)/T + s \\ \vdots \\ (g_N - k)/T + s \\ u/T \\ -1/T \end{bmatrix}  =
    \begin{bmatrix} -s_1 + r_1 + s \\ \vdots \\ - s_N + r_N + s \\ 0 \\ -\frac{\rho c_v}{T} \end{bmatrix}, \
    \bigg(\frac{\partial \mathbf{u}}{\partial Z}\bigg)^{\tr}
         \frac{\partial \mathbf{v}}{\partial Z}
        = 
        \begin{bmatrix} 
            r_1/\rho_1  &        &     0      &    0   &    0   \\
                        & \ddots &            & \vdots & \vdots \\
            0           &        & r_N/\rho_N &    0   &    0   \\
            0           & \hdots &     0      & \rho/T &    0 \\
            0           & \hdots &     0      &    0   &  \rho c_v /T^2 
        \end{bmatrix}.    
\end{equation*}
Therefore:
\begin{equation}\label{eq:H2}
    \mathbf{H}_{II} = f^{'} \begin{bmatrix} 
            r_1/\rho_1  &        &     0      &    0   &    0   \\
                        & \ddots &            & \vdots & \vdots \\
            0           &        & r_N/\rho_N &    0   &    0   \\
            0           & \hdots &     0      & \rho/T &    0 \\
            0           & \hdots &     0      &    0   &  \rho c_v /T^2 
        \end{bmatrix} - \frac{f^{''}}{\rho} \begin{bmatrix} R_1 \\ \vdots \\ R_N \\ 0 \\ -\frac{\rho c_v}{T} \end{bmatrix} \begin{bmatrix} R_1 & \hdots & R_N & 0 & -\frac{\rho c_v}{T} \end{bmatrix}, 
\end{equation}
where $R_i = - s_i + r_i + s$. We recover Harten's conditions in the compressible Euler case. At this point, we immediately note that if $ f^{'} > 0, \ f^{''} < 0$ then $\mathbf{H}_{II}$ is positive definite (as the sum of a positive definite matrix and a positive semi-definite matrix). Therefore a minimum entropy principle for the mixture's specific entropy holds. \\
\indent Continuing on the characterization of convexity, $\mathbf{H}_{II}$ writes:
\begin{equation*}
    \mathbf{H}_{II} = \frac{f^{'}}{\rho} \begin{bmatrix} 
            r_1/Y_1  &        &     0      &    0   &    0   \\
                        & \ddots &            & \vdots & \vdots \\
            0           &        & r_N/Y_N &    0   &    0   \\
            0           & \hdots &     0      & \rho^2/T &    0 \\
            0           & \hdots &     0      &    0   &  \rho^2 c_v /T^2 
        \end{bmatrix} - \frac{f^{''}}{\rho} \begin{bmatrix}
        R_1^2                   &        &         R_1 R_N          &    0   & -\frac{\rho c_v}{T} R_1    \\
                                & \ddots &                          & \vdots &         \vdots             \\
        R_1 R_N                 &        &         R_N^{2}          &    0   & -\frac{\rho c_v}{T} R_N    \\
        0                       & \hdots &            0             &    0   &            0               \\
        -\frac{\rho c_v}{T} R_1 & \hdots & -\frac{\rho c_v}{T} R_N  &    0   & \frac{\rho^2 c_v^{2} }{T^2}
        \end{bmatrix}.
\end{equation*}
Let $\overline{r}_i = r_i/Y_i$ and $\eta = f^{'} - c_{v}f^{''}$, for two species we have:
\begin{equation*}
    \mathbf{H}_{II} = \frac{1}{\rho}\begin{bmatrix} 
        f^{'} \overline{r}_1 - f^{''}R_1^2 & -R_1 R_2 f^{''} & 0 & \frac{\rho c_v}{T} R_1 f^{''}\\
        -R_1 R_2 f^{''} &f^{'} \overline{r}_2 - f^{''} R_2^2  &  0 & \frac{\rho c_v}{T} R_2 f^{''} \\
        0 & 0 & \rho^2 f^{'}/T & 0 \\
        \frac{\rho c_v}{T} R_1 f^{''} & \frac{\rho c_v}{T} R_2 f^{''} & 0 & \frac{\rho^2 c_v}{T^2} \eta
        \end{bmatrix}
\end{equation*}
The determinants of the first three major blocks of $\mathbf{H}$ are:
\begin{equation}\label{eq:major3H}
    H_{11} = \overline{r}_1 \bigg( f^{'}  - f^{''} \frac{R_1^2}{\overline{r}_1}\bigg), \ H_{22} = \overline{r}_1 \overline{r}_2  f^{'} \bigg( f^{'} -  f^{''} \bigg( \frac{R_1^2}{\overline{r}_1}  + \frac{R_2^2}{\overline{r}_2} \bigg) \bigg), \ H_{33} = \frac{\rho^2 f^{'}}{T} H_{22}.
\end{equation}
Last:
\begin{align*}
    det(\rho \mathbf{H}_{II}) =& \ \frac{\rho^2 f^{'}}{T} 
        \begin{vmatrix} 
        f^{'} \overline{r}_1 - f^{''}R_1^2 & - R_1 R_2 f^{''}                & \frac{\rho c_v}{T} R_1 f^{''} \\
        -R_1 R_2 f^{''}               &  f^{'} \overline{r}_2 - f^{''} R_2^2 & \frac{\rho c_v}{T} R_2 f^{''} \\
        \frac{\rho c_v}{T} R_1 f^{''} & \frac{\rho c_v}{T} R_2 f^{''} & \frac{\rho^2 c_v}{T^2}\eta
        \end{vmatrix} \\
                    =& \ \frac{\rho^4 c_v f^{'}}{T^3} 
        \begin{vmatrix} 
        f^{'} \overline{r}_1 - f^{''}R_1^2 & - R_1 R_2 f^{''} & R_1 f^{''}\\
        - R_1 R_2 f^{''} &f^{'} \overline{r}_2 - f^{''} R_2^2 & R_2 f^{''}\\
        c_v R_1 f^{''} & c_v R_2 f^{''}& \eta
        \end{vmatrix} \\
                    =& \ \frac{\rho^4 c_v f^{'}}{T^3} \bigg( \eta H_{22} - 
                           c_v R_2^2 f^{''} \begin{vmatrix} 
                                    f^{'} \overline{r}_1 - f^{''}R_1^2 & - R_1 f^{''}\\
                                    R_1  & 1
                                \end{vmatrix} -
                            c_v R_1^2 f^{''} \begin{vmatrix} 
                                    f^{'} \overline{r}_2 - f^{''}R_2^2 & - R_2 f^{''}\\
                                    R_2 & 1
                                \end{vmatrix}    
                            \bigg) \\
                    =& \ \frac{\rho^4 c_v f^{'}}{T^3} \bigg( \eta H_{22} - 
                           c_v f^{''} f^{'} ( R_2^2  \overline{r}_1 + R_1^2 \overline{r}_2 \big)  
                            \bigg) \\
                    =& \ \frac{\rho^4 c_v (f^{'})^2}{T^3} \overline{r}_1 \overline{r}_2 \bigg( \eta   f^{'} - (\eta + c_{v}) f^{''} \bigg( \frac{R_1^2}{\overline{r}_1} + \frac{R_2^2}{\overline{r}_2} \bigg) \bigg).
\end{align*}
We obtain conditions on $f$ involving terms of the form $f^{'} - \alpha f^{''}$, but unlike in the Euler case, $\alpha$ is not a constant. In section 4.1, the simple structure of the mapped Hessian $\mathbf{H}_{I}$, given by equation (\ref{eq:H1}), allowed us to easily derive the necessary and sufficient conditions (\ref{eq:Harten_MC}) for convexity for an arbitrary number of species. Nevertheless, we were not able to conclude on a minimum entropy principle on the specific entropy of each species. Here, the mapped Hessian $\mathbf{H}_{II}$, given by equation  (\ref{eq:H2}), is mostly dense, which complicates the task of establishing convexity conditions for an arbitrary number of species. However, we know from equation (\ref{eq:H2}) that $f^{'} > 0$ and $f^{''} < 0$ are sufficient conditions for admissibility, independently of the number of species, which is enough to conclude on a minimum entropy principle on the mixture's specific entropy.

\section{Numerical schemes satisfying a minimum entropy principle}\label{sec:numerics}
In this section, we review schemes which, by virtue of satisfying all entropy inequalities under some assumptions, satisfy a minimum entropy principle for the compressible multicomponent Euler equations. \\
\indent We only discuss first-order schemes in one dimension. Extensions to high-order and multiple dimensions (including unstructured grids) can be found in \cite{Zhang, Ihme, Guermond_IDP, Guermond_IDP2, Guermond_IDP3}. These schemes are typically constructed as composite \textit{convex} combinations of one-dimensional first-order updates. Since entropies are convex functions, any entropy inequality satisfied by the baseline one-dimensional first-order update will be satisfied by the whole scheme as well. 
\subsection{Godunov-type schemes \cite{UPS}}
Let $\mathbf{w}(x/t; \mathbf{u}_L, \mathbf{u}_R)$ be the solution of the Riemann problem:
\begin{equation}\label{eq:Riemann_pb}
    \partial_t \mathbf{u}  + \partial_x \mathbf{f} = 0, \ 
        \mathbf{u}(x, 0) = \left\{
                \begin{array}{ll}
                  \mathbf{u}_L, \ \ x < 0, \\
                  \mathbf{u}_R, \ \ x > 0, \\
                \end{array}
              \right.
\end{equation}
where $\mathbf{u}_L$ and $\mathbf{u}_R$ are constant states. Let $a_L$ and $a_R$ be the smallest and largest signal velocities. Then $\mathbf{w}$ satisfies:
\begin{equation}\label{eq:Riemann}
    \mathbf{w}(x/t; \mathbf{u}_L, \mathbf{u}_R) = \left\{
                \begin{array}{ll}
                  \mathbf{u}_L, \ \ x/t \leq a_L \\
                  \mathbf{u}_R, \ \ x/t \geq a_R \\
                \end{array}
              \right. 
\end{equation}
In the Godunov scheme \cite{Godunov}, each discontinuity in the discrete field $\mathbf{u}_i^n$ gives rise to a local Riemann problem (\ref{eq:Riemann}). If $\lambda |a_{max}| < 1/2$, where $a_{max}$ is the largest signal speed in the domain, then there is no interaction between neighboring Riemann problems and the exact solution $\mathbf{w}_{n+1}(x)$ at the next time instant writes:
\begin{equation*}
    \mathbf{w}_{n+1}(x) = \mathbf{w}((x - x_{i+\frac{1}{2}} )/\Delta t ; \mathbf{u}_{i}^{n}, \mathbf{u}_{i+1}^{n}), \ \mbox{for} \ |x - x_{i+\frac{1}{2}}| \leq \Delta x / 2,
\end{equation*}
where $x_{i+\frac{1}{2}}$ is the position of the interface between cells $i$ and $i+1$. The Godunov scheme is obtained by averaging $\mathbf{w}_{n+1}$ in each cell:
\begin{align*}
    \mathbf{u}_i^{n+1} =& \ \frac{1}{\Delta x} \int_{x_{i-\frac{1}{2}}}^{x_{i+\frac{1}{2}}} \mathbf{w}_{n+1}(x) \ dx \\ 
                       =& \ \frac{1}{\Delta x} \int_{0}^{ \Delta x/2} \mathbf{w}(x/\Delta t; \mathbf{u}^{n}_{i-1}, \mathbf{u}^{n}_{i} ) \ dx +
                            \frac{1}{\Delta x} \int_{-\Delta x/2}^{0} \mathbf{w}(x/\Delta t; \mathbf{u}^{n}_{i}, \mathbf{u}^{n}_{i+1} ) \ dx.
\end{align*}
This update can be rewritten in conservative form: 
\begin{equation*}
    \mathbf{u}_i^{n+1} = \mathbf{u}_i^{n} - \lambda \big( \mathbf{f}(\mathbf{\hat{w}}_{i+\frac{1}{2}}) - 
                                      \mathbf{f}(\mathbf{\hat{w}}_{i-\frac{1}{2}})   \big), \ \mathbf{\hat{w}}_{i+\frac{1}{2}} = \mathbf{w}(0; \mathbf{u}_i^{n}, \mathbf{u}_{i+1}^n),
\end{equation*}
with $\lambda = \Delta t/\Delta x$. An important assumption from there \cite{UPS, Guermond_IDP} is that the exact Riemann solution is an entropy solution. This implies, for all entropies:
\begin{equation*}
    \frac{1}{\Delta x}\int_{x_{i-\frac{1}{2}}}^{x_{i+\frac{1}{2}}} U(\mathbf{w}_{n+1}(x)) \ dx \leq U(\mathbf{u}_i^n) - \lambda \big( F (\mathbf{\hat{w}}_{i+\frac{1}{2}}) - F (\mathbf{\hat{w}}_{i-\frac{1}{2}}) \big). 
\end{equation*}
With Jensen's inequality:
\begin{equation*}
    U\bigg(\frac{1}{\Delta x} \int_{x_{i-\frac{1}{2}}}^{x_{i+\frac{1}{2}}} \mathbf{w}_{n+1}(x) \ dx \bigg) \leq \frac{1}{\Delta x} \int_{x_{i-\frac{1}{2}}}^{x_{i+\frac{1}{2}}} U(\mathbf{w}_{n+1}(x)) \ dx,
\end{equation*}
it follows that the Godunov scheme satisfies:
\begin{equation}\label{eq:Godunov_ES}
    U(\mathbf{u}_{i}^{n+1}) \leq U(\mathbf{u}_i^n) - \lambda \big( F (\mathbf{\hat{w}}_{i+\frac{1}{2}}) - F (\mathbf{\hat{w}}_{i-\frac{1}{2}}) \big). 
\end{equation}
This shows that the Godunov scheme inherits, by construction, all the entropy inequalities that the exact Riemann solution satisfies. This result also applies to schemes based on approximate Riemann solutions provided that they remain consistent with the integral forms of the conservation law and the entropy inequality (see Theorem 3.1 in \cite{UPS}). The bottom line is that \textit{full knowledge of the Riemann solution is not necessary}. For instance, the HLL scheme \cite{UPS} qualifies if the maximum right and left wave speeds are correctly estimated (from above). \\
\indent The Godunov scheme satisfies a sharper version of (\ref{eq:min_entropy}). Using (\ref{eq:Godunov_ES}) with $U = - \rho f_0(s)$ and $ s_0 = \min [s(\mathbf{u}_{i-1}^{n}), \ s(\mathbf{u}_{i}^{n}), \ s(\mathbf{u}_{i+1}^{n})]$, and the fact that the exact solution $\mathbf{w}$ is an entropy solution satisfying (\ref{eq:min_entropy_OP}), it follows that the Godunov scheme satisfies:
\begin{equation}
    s(\mathbf{u}_i^{n+1}) \geq \min [s(\mathbf{u}_{i-1}^{n}), \ s(\mathbf{u}_{i}^{n}), \ s(\mathbf{u}_{i+1}^{n})],
\end{equation}
\indent For the compressible Euler equations, procedures for calculating the exact solution (see Toro \cite{Toro}) and estimating the maximum wave speed (see Guermond \& Popov \cite{Guermond_speed}) are available and can be extended to the multicomponent case (a follow-up to \cite{Guermond_speed} is proposed by Frolov in \cite{Frolov}, section 4.5). \\
\indent It is unclear whether the assumption that the exact Riemann solution satisfies \textit{all} entropy inequalities is valid. To the best of the authors' knowledge, there is no proof that Harten's entropies \cite{Harten} are the only entropies of the compressible Euler equations. The same can be said about the entropies that we explored in section \ref{sec:entropy_functions} for the multicomponent case. This precludes a direct proof where entropy inequalities are evaluated for the exact Riemann solution. Another way of proving this would be to show that the exact Riemann solution can be written as a limit solution to the regularized system (\ref{eq:PDE_visc_1D}) or any other sytem which implies all entropy inequalities. As far as the minimum entropy principle is concerned, showing that the exact Riemann solution satisfies all entropy inequalities associated with Harten's family or with the convolution entropies of section \ref{sec:minS_proof} would be enough. 
\subsection{The Lax-Friedrichs scheme}
The Lax-Friedrichs (LxF) scheme writes:
\begin{equation*}
    \mathbf{u}_{i}^{n+1} = \frac{\mathbf{u}_{i-1}^{n} + \mathbf{u}_{i+1}^{n}}{2} + \frac{\lambda}{2}\big(\mathbf{f}(\mathbf{u}_{i-1}^{n}) - \mathbf{f}(\mathbf{u}_{i+1}^{n})\big).
\end{equation*}
\indent Harten (private communication in \cite{Tadmor4}, section 4) observed that if the time step is small enough, the LxF scheme coincides with the Godunov scheme over a staggered grid. The solution thus inherits the entropy inequalities that the Riemann solution satisfies:
\begin{equation}\label{eq:LxF_ES}
    U(\mathbf{u}_{i}^{n+1}) \leq \frac{U(\mathbf{u}_{i-1}^{n}) + U(\mathbf{u}_{i+1}^{n})}{2} + \frac{\lambda}{2}\big(F(\mathbf{u}_{i-1}^{n}) - F(\mathbf{u}_{i+1}^{n})\big).
\end{equation}
As in section \ref{sec:minS_proof}, it is easy to show that inequality (\ref{eq:LxF_ES}) with $U = - \rho f_0(s)$ and $ s_0 = \min [s(\mathbf{u}_{i-1}^{n}), \ s(\mathbf{u}_{i+1}^{n})]$ leads to a minimum entropy principle:
\begin{equation}
    s(\mathbf{u}_i^{n+1}) \geq \min [s(\mathbf{u}_{i-1}^{n}), \ s(\mathbf{u}_{i+1}^{n})],
\end{equation}
that is sharper than (\ref{eq:min_entropy}). \\
\indent On the other hand, Lax \cite{Lax} proved, without invoking Riemann solutions, that the LxF scheme can be made to satisfy (\ref{eq:LxF_ES}) for \textit{any} given entropy pair. We recall his proof here, as it will help us address a point brought up during the review process. \\
\indent Denote $\mathfrak{u} = \mathbf{u}_{i}^{n+1}$, $\mathfrak{v} = \mathbf{u}_{i-1}^{n}$ and $\mathfrak{w} = \mathbf{u}_{i+1}^{n}$. The LxF scheme writes: 
\begin{equation*}
    \mathfrak{u}(\mathfrak{v}, \mathfrak{w}) = \frac{\mathfrak{v} + \mathfrak{w}}{2} + \frac{\lambda}{2}(\mathbf{f}(\mathfrak{v}) - \mathbf{f}(\mathfrak{w})), 
\end{equation*}
and the entropy inequality (\ref{eq:LxF_ES}) can be studied by looking at the sign of the difference function:
\begin{equation*}
    \Delta\mathcal{S}(\mathfrak{v}, \mathfrak{w}) = \frac{U(\mathfrak{v}) + U(\mathfrak{w})}{2} + \frac{\lambda}{2}(F(\mathfrak{v}) - F(\mathfrak{w})) - U(\mathfrak{u}). 
\end{equation*}
Lax \cite{Lax} used a homotopy approach. Let $s \in [0 \ 1]$, and define:
\begin{equation*}
    \overline{\mathfrak{v}}(s) = s \mathfrak{v} + (1 - s) \mathfrak{w}, \ \overline{\mathfrak{u}}(s) = \mathfrak{u}(\overline{\mathfrak{v}}(s), \mathfrak{w}).
\end{equation*}
Since $\overline{\mathfrak{v}}(1) = \mathfrak{v}, \ \overline{\mathfrak{v}}(0) = \mathfrak{w}$, and $\Delta\mathcal{S}(\mathfrak{w}, \mathfrak{w}) = 0$, the fundamental theorem of calculus gives:
\begin{equation}\label{eq:Lax_proof_1}
    \Delta\mathcal{S}(\mathfrak{v}, \mathfrak{w}) = \Delta\mathcal{S}(\overline{\mathfrak{v}}(1), \mathfrak{w}) - \Delta\mathcal{S}(\overline{\mathfrak{v}}(0), \mathfrak{w}) = \int_{0}^1 \frac{d}{ds} \bigg( \Delta\mathcal{S}(\overline{\mathfrak{v}}(s), \mathfrak{w}) \bigg) \ ds.
\end{equation}
$\overline{\mathfrak{u}}$ and $\overline{\mathfrak{v}}$ satisfy:
\begin{equation*}
    \frac{d \overline{\mathfrak{v}}}{ds} = \mathfrak{v} - \mathfrak{w}, \ \frac{d \overline{\mathfrak{u}}}{ds} = \frac{\mathfrak{v} - \mathfrak{w}}{2} + \frac{\lambda}{2}A(\overline{\mathfrak{v}})(\mathfrak{v} - \mathfrak{w}) = \frac{1}{2} \bigg( I + \lambda A(\overline{\mathfrak{v}})  \bigg) (\mathfrak{v} - \mathfrak{w}),
\end{equation*}
where $A$ is the flux Jacobian. Using chain rules and the constitutive relation (\ref{eq:entropy_cons}), the integrand in equation (\ref{eq:Lax_proof_1}) writes:
\begin{equation*}
    \frac{d}{ds} \bigg( \Delta\mathcal{S}(\overline{\mathfrak{v}}(s), \mathfrak{w}) \bigg) = \frac{1}{2}  \bigg( \frac{dU}{d\mathbf{u}}(\overline{\mathfrak{v}}) - \frac{dU}{d\mathbf{u}}(\overline{\mathfrak{u}}) \bigg) \bigg(I + \lambda A(\overline{\mathfrak{v}}) \bigg) (\mathfrak{v} - \mathfrak{w}).
\end{equation*}
Again, let $r \in [0 \ 1]$, and define:
\begin{equation*}
    \overline{\mathfrak{w}}(r,s) = r \overline{\mathfrak{v}}(s) + (1 - r) \mathfrak{w} = r s \mathfrak{v} + (1 - r s) \mathfrak{w}, \ \overline{\overline{\mathfrak{u}}}(r,s) = \mathfrak{u}(\overline{\mathfrak{v}}(s), \overline{\mathfrak{w}}(r)).
\end{equation*}
Since $\overline{\overline{\mathfrak{u}}}(1,s) = \overline{\mathfrak{v}}(s), \ \overline{\overline{\mathfrak{u}}}(0,s) = \overline{\mathfrak{u}}(s)$, the fundamental theorem of calculus gives:
\begin{equation}\label{eq:Lax_proof2}
    \frac{dU}{d\mathbf{u}}(\overline{\mathfrak{v}}) - \frac{dU}{d\mathbf{u}}(\overline{\mathfrak{u}}) = \int_{0}^{1} \frac{d}{dr} \bigg( \frac{dU}{d\mathbf{u}}(\overline{\overline{\mathfrak{u}}})  \bigg) \ dr = \int_{0}^1 \bigg( \frac{d\overline{\overline{\mathfrak{u}}}}{dr} \bigg)^T G(\overline{\overline{\mathfrak{u}}})  \ dr,
\end{equation}
where $G$ is the entropy Hessian. With:
\begin{equation*}
    \frac{d\overline{\overline{\mathfrak{u}}}}{dr} = \frac{s}{2} \bigg( I - \lambda A(\overline{\mathfrak{w}})\bigg) (\mathfrak{v} - \mathfrak{w})
\end{equation*}
and equations (\ref{eq:Lax_proof_1}) - (\ref{eq:Lax_proof2}), the difference function $\Delta \mathcal{S}$ can finally be rewritten as:
\begin{align*}
    \Delta\mathcal{S}(\mathfrak{v}, \mathfrak{w}) =& \ \int_{0}^1 \int_{0}^1 \frac{s}{4}  \bigg( \big( I - \lambda A(\overline{\mathfrak{w}})\big) (\mathfrak{v} - \mathfrak{w}) \bigg)^{T} G(\overline{\overline{\mathfrak{u}}}) \bigg( \big(I + \lambda A(\overline{\mathfrak{v}}) \big)  (\mathfrak{v} - \mathfrak{w}) \bigg) \ ds dr. \\
                                                  =& \ \langle \mathfrak{z}, \ \mathfrak{z} \rangle_G  
    - \lambda ( \langle A(\overline{\mathfrak{w}}) \mathfrak{z}, \ \mathfrak{z} \rangle_G  +
                \langle \mathfrak{z} , \ A(\overline{\mathfrak{v}}) \mathfrak{z} \rangle_G ) 
    - \lambda^2 \langle A(\overline{\mathfrak{w}}) \mathfrak{z}, \ A(\overline{\mathfrak{v}}) \mathfrak{z} \rangle_G.   
\end{align*}
where $\mathfrak{z} = (\mathfrak{v} - \mathfrak{w})$ and $\langle \ , \  \rangle_G$ is the inner product defined by:
\begin{equation*}
    \langle \mathfrak{a}, \ \mathfrak{b} \rangle_G = \int_{0}^1 \int_0^{1} \frac{s}{4} \mathfrak{a}^T G(\overline{\overline{u}}) \mathfrak{b} \ dsdr.
\end{equation*}
Since $G$ is symmetric positive definite, $\langle \mathfrak{z}, \ \mathfrak{z} \rangle_G > 0$ and one can expect the entropy inequality (\ref{eq:LxF_ES}) to be met if $\lambda$ is small enough. Within the vector space spanned by $(r, \ s)$, let $c$ be the maximum matrix norm of $A$, $m$ be the minimum eigenvalue of $G$ and $M$ be the maximum eigenvalue of $G$. Then, for $||\mathfrak{v}|| \neq 0$, if $\lambda$ satisfies:
\begin{equation}\label{eq:CFLL}
    m - 2 c \lambda M - c^2 \lambda^2 M > 0 \ \Leftrightarrow \ \lambda c < \sqrt{1 + (m/M)} - 1. 
\end{equation}
then the inequality (\ref{eq:LxF_ES}) is met. Since $U$ is strictly convex, $(m/M) > 0$ and the right-hand side of (\ref{eq:CFLL}) is strictly positive. In other words, for \textit{any} entropy $U$, there will always exist a time step small enough such that the condition (\ref{eq:CFLL}) is met. \\
\indent While Lax's proof does not invoke Riemann solutions, it does not completely support the statement \cite{Tadmor2} that the LxF scheme can be made to satisfy \textit{all entropy inequalities}. The factor $m/M$ in (\ref{eq:CFLL}) is strictly positive, but also depends on the entropy at hand. The fact that we do not know all the entropies of a hyperbolic system in general leaves open the possibility that $m/M$ can be arbitrarily small. One needs to show that there exists a strictly positive and entropy-independent lower bound $K$ on $m/M$, so that under the condition:
\begin{equation}\label{eq:CFLK}
     \lambda c < \sqrt{1 + K} - 1
\end{equation}
the LxF scheme will effectively satisfy all entropy inequalities. As far as the minimum entropy principle is concerned however, we recalled in section \ref{sec:minS_proof} that not all entropy inequalities need to be satisfied.

\section{Conclusions}
\hspace{0.4 cm} We proved a minimum entropy principle for entropy solutions to the multicomponent compressible Euler equations, extending Tadmor's result \cite{Tadmor2}. The proof was carried out in one dimension but easily follows in two and three dimensions (the characterization of the two families in section \ref{sec:entropy_functions} is independent of the number of dimensions). This principle was proven for the mixture's specific entropy only. It would be interesting to establish whether this also holds for the specific entropy of each species. We assumed a mixture of thermally perfect gases governed by an ideal gas law. The methodology outlined here and in the work of Harten \textit{et al.} \cite{SuperHarten}, which extended Harten's characterization \cite{Harten} to gases with an arbitrary equation of state, should provide helpful guidelines for those interested in taking this result farther. \\ 
\indent While numerical schemes consistent with the entropy condition (\ref{eq:entropy_cond}) for a given pair $(U,  F)$ can be constructed \cite{Tadmor} (for the compressible multicomponent Euler equations, Gouasmi \textit{et al.} \cite{Gouasmi} constructed one such scheme for the pair $(-\rho s, -\rho u s)$), designing numerical schemes which lead to discrete entropy solutions is more challenging. A common trait of such schemes \cite{UPS, Zhang, Guermond_IDP, Ihme} is that they take root in the notion of a Riemann problem and the existence of solutions satisfying all entropy inequalities. \\
\indent While the minimum entropy principle is only a property of entropy solutions, it provides valuable information about the \textit{local} behavior of the physical solution. Limiting procedures for high-order schemes have been designed around this property \cite{Zhang, Ihme, Guermond_IDP2, Guermond_IDP3} for the Euler equations and may henceforth prove useful in multicomponent flow simulations. \\
\indent Finally, we emphasize that the present work is not meant to provide a comprehensive review of the symmetrizability of the multicomponent system. We refer the interested reader to Giovangigli \& Matuszewski \cite{Giovangigli_Mat} for instance. The investigation of entropy functions carried out in section \ref{sec:entropy_functions} was driven by the prospect of proving a minimum entropy principle. Harten's pioneering work \cite{Harten} had broader motivations. 

\ifx
\section*{Acknowledgments}
\hspace{0.4 cm} Ayoub Gouasmi and Karthik Duraisamy were funded by the AFOSR through grant number FA9550-16-1-030 (Tech. monitor: Fariba Fahroo). Eitan Tadmor was supported in part by NSF grants DMS16-13911, RNMS11-07444 (KI-Net) and ONR grant N00014-1812465. \\
\indent Ayoub Gouasmi is grateful to Jean-Luc Guermond, Bojan Popov and Ignacio Tomas for spirited discussions on their work \cite{Guermond_visc, Guermond_IDP, Guermond_IDP2} and for bringing the work of Delchini \textit{et al.} \cite{Delchini_1, Delchini_2} and Harten \textit{et al.} \cite{SuperHarten} to his attention. 
\fi

\end{document}